# DISCUSSION OF "ANALYSIS OF VARIANCE—WHY IT IS MORE IMPORTANT THAN EVER" BY A. GELMAN

By Peter McCullagh[1]

*University of Chicago*

Factorial models and analysis of variance have been a central pillar of statistical thinking and practice for at least 70 years, so the opportunity to revisit these topics to look for signs of new life is welcome. Given its lengthy history, ANOVA is an unpromising area for new development, but recent advances in computational techniques have altered the outlook regarding spatial factors, and could possibly affect what is done in nonspatial applications of factorial designs. In common with a large number of statisticians, Gelman regards analysis of variance as an algorithm or procedure with a well-defined sequence of computational steps to be performed in fixed sequence. The paper emphasizes tactics, how to "perform ANOVA," how to "set up an ANOVA," how to compute "the correct ANOVA," what software to use and how to use it to best effect. The "solution to the ANOVA problem" proffered in Section 3.2 emphasizes once again, how to do it in the modern hierarchical style. Were it not for the recommendation favoring shrinkage, one might have expected a more accurate descriptive title such as the Joy of ANOVA.

I admire the breezy style, the fresh approach and the raw enthusiasm of this paper. It contains perhaps three points with which I agree, namely the importance of ANOVA, the usefulness of thinking in terms of variance components and a passage in Section 3.4 on near-zero estimated variance components. How we could agree on these points and disagree on nearly everything else takes a good deal of explanation. My own philosophy is that it is usually best to begin with the question or questions, and to tailor the analyses to address those questions. Generally speaking, one expects the answer to depend on the question, and it is unreasonable to ask that the analysis, or even a major part of it, should be the same for all questions asked. In my experience, routine statistical questions are less common than questionable statistical routines, so I am loath to make pronouncements about what is and what is not relevant in applied statistics. In one of his least convincing passages Gelman argues that the new methodology does the right thing automatically, even for complicated designs. I am inclined to regard this claim either as a regrettable rhetorical flourish, or as a self-fulfilling statement *defining* the class of designs and factors with which the paper is concerned. In the latter case, there is little left to discuss, except to protest that large segments of analysis of variance and factorial design have been overlooked.

[1]Supported in part by NSF Grant 03-05009.

The phrase "random coefficient model" or "varying coefficient model" is one that ought to trigger alarm bells. If $x$ is temperature in °C and $x'$ is temperature in °F, the linear models

$$\beta_0 + \beta_1 x \quad \text{and} \quad \beta'_0 + \beta'_1 x'$$

are equivalent in the sense that they determine the same subspace and thus the same set of probability distributions. Consider now the model in which $\beta_1 \sim N(\bar{\beta}_1, \sigma_1^2)$ and the corresponding one in which $\beta'_1 \sim N(\bar{\beta}'_1, \tau_1^2)$. On the observation space, the implied marginal covariances are

$$\sigma^2 I_n + \sigma_1^2(xx^\top) \quad \text{and} \quad \sigma^2 I_n + \tau_1^2(x'x'^\top),$$

two linear combinations of matrices spanning different spaces. In other words, these random-coefficient formulations do not determine the same set of distributions. It is only in very special circumstances that a random-effects model constructed in this way makes much sense. Making sense is a property that is intuitively obvious: mathematically it means that the model is a group homomorphism or representation.

Gelman's paper is concerned almost exclusively with simple factorial designs in which the factor effects are plausibly regarded as exchangeable. A batch is not a set of regression coefficients as suggested in Section 3.2, but a set of *effects*, one effect for each factor level, and one set or batch for each factor or interaction. The preceding paragraph shows why the distinction between coefficient and effect matters. If batch were synonymous with subset, the new term would be redundant, so it appears that the effects in a batch are meant to be random. In Section 6, a batch of effects is defined as a set of random variables, which are then treated as exchangeable without comment, as if no other option exists. The grouping by batches is determined by factor levels, which is automatic for simple factorial designs, nested or crossed. However, this is not necessarily the case for more general factorial structures such as arise in fertility studies [Cox and Snell (1981), pages 58–62], tournament models [Joe (1990)], origin-destination designs [Stewart (1948)], import-export models, citation studies [Stigler (1994)] or plant breeding designs in which the same factor occurs twice.

In virtually all of the literature on factorial design and analysis of variance, effects are either fixed or random. No other types are tolerated, and all random effects are independent and identically distributed, as in Section 6 of the present paper. This regrettable instance of linguistic imperialism makes it difficult to find a satisfactory term for random effects in which the components are random but not independent. Clarity of language is important, and in this instance the jargon has developed in such a way that it has become a major obstacle to communication. My own preference is to address matters of terminology, such as treatment and block factors, fixed and random effects, and so on, by what they imply in a statistical model, as described in the next two paragraphs. The alternative to these definitions is the linguistic quagmire so well documented by Gelman in Section 6.

A treatment factor or classification factor $A$ is a list such that $A(i)$ is the level of factor $A$ on unit $i$. Usually, the set of levels is finite, and the information may then be coded in an indicator matrix $X = X(A)$, one column for each level. By contrast, a block factor $E$ is an equivalence relation on the units such that $E_{ij} = 1$ if units $i, j$ are in the same block, and zero otherwise. A treatment or classification factor may be converted into a block factor by the forgetful transformation $E = XX^\top$ in which the names of the factor levels are lost. A block factor cannot be transformed into a treatment factor because the blocks are unlabelled. A factor may occur in a linear model in several ways, the most common of which are additively in the mean and additively in the covariance

$$(1) \qquad Y \sim N(X\beta, \sigma^2 I_n) \quad \text{or} \quad Y \sim N(\mathbf{1}\mu, \sigma^2 I_n + \sigma_b^2 E).$$

Traditionally, the terms "fixed-effects model" and "random-effects model" are used here, but this terminology is not to be encouraged because it perpetuates the myth that random effects are necessarily independent and identically distributed. Note that $I_n$ is the equivalence relation corresponding to units, and $\sigma^2 I_n$, the variance of the exchangeable random unit effects, is included in both models.

Suppose now that two factors $A, B$ are defined on the same set of units, and that these factors are crossed, $A.B$ denoting the list of ordered pairs. The corresponding block factors may be denoted by $E_A$, $E_B$ and $E_{AB}$. Two factors may occur in a linear model in several ways, the conventional factorial models for the mean being denoted by

$$1, \quad A, \quad B, \quad A+B, \quad A.B,$$

with a similar list of linear block-factor models for the covariances

$$I, \quad I + E_{AB}, \quad I + E_A, \quad I + E_B,$$

$$I + E_A + E_B, \quad I + E_A + E_B + E_{AB}.$$

Here $A + B$ denotes the vector space of additive functions on the factor levels, whereas $I + E_A + E_B$ denotes the set of nonnegative combinations of three matrices in which the coefficients are called variance components. If a factor occurs in the model for the mean, the associated variance component is not identifiable. For example, if the model for the mean includes $A.B$, a so-called nonrandom interaction, the variance components associated with $E_A, E_B, E_{AB}$ are not identifiable. However, if the variance model includes the interaction $E_{AB}$, the additive model $A + B$ for the mean is ordinarily identifiable. These are mathematical statements concerning the underlying linear algebra. Philosophical pronouncements such as "if one main effect is random, the interaction is also random" have no place in the discussion.

The subspace $A \subset \mathcal{R}^n$ determined by a factor is of a very special type: it is also a ring, closed under functional multiplication, with **1** as identity element. A factorial model is also a special type of vector subspace of functions on the units. Each is a representation of the product symmetric group in the tensor product space that is also closed under deletion of levels [McCullagh (2000)]. Each of the variance-component models listed above is also a representation in the same sense, but one in which the subspace consists of certain symmetric functions on ordered pairs of units, that is, symmetric matrices. Specifically, each exchangeable variance-components model is a *trivial* representation in the space of symmetric matrices that is closed under deletion of levels. By contrast, a Taguchi-type model in which the variance depends on one or more factor levels is a representation, but not a trivial representation. This may not be a helpful statement for most student audiences, but it does serve to emphasize the point that factorial subspaces are determined by groups and representations. ANOVA decomposition requires one further ingredient in the form of an inner product on the observation space.

If the term "classical linear regression model" implies independence of components, as Gelman's usage in Section 3.3 suggests, then most of the factorial models described above are not classical. On the other hand, they have been a part of the literature in biometry and agricultural field trials for at least 70 years, so they are not lacking in venerability. For clarity of expression, the term "neoclassical" is used here to include models of the above type, linear in the mean and linear in the covariance. The prefix "neo-" refers to more recent versions, including certain spatial models, spline-smoothing models and Taguchi-type industrial applications in which the primary effect of so-called noise factors [Wu and Hamada (2000)] is on variability. A pure variance-components model is one in which the model for the mean is trivial, that is, the constant functions only. The simplest neoclassical procedure for estimation and prediction is first to compute the variance components using residual maximum likelihood, then to compute regression coefficients by weighted least squares, and then to compute predicted values and related summary statistics. For prediction to be possible, the model must be a family of processes.

The main thrust of Gelman's paper as I understand it is to argue that ANOVA should be performed and interpreted in the context of an additive variance-components model rather than an additive factorial model for the mean. This is the special neoclassical model in which the subspace for the mean is the one-dimensional vector space of constant functions, and all effects and interactions are included as block factors in the covariance function. A joint prior distribution on the variance components avoids the discontinuity associated with near-zero estimated variance components. Individual treatment effects do not occur as parameters in this model, but they may be estimated by prediction, that is, by computing the conditional mean for a new unit having a given factor level, or the difference between conditional means for two such units. When the factor levels are numerous or nonspecific [Cox (1984)], or ephemeral or faceless [Tukey (1974)], this approach is uncontroversial, and indeed, strongly recommended. However, numerous examples exist in which one or more factors have levels that are not of this type, where inference for a specific treatment contrast or a specific classification contrast is the primary purpose of the experiment. Exchangeability is simply one of many modeling options, sensible in many cases, debatable in others, and irrelevant for the remainder. To my mind, Gelman has failed to make a convincing case that additive models for the mean should be abandoned in favor of a scheme that "automatically gets it right" by assuming that every factor has levels whose effects are exchangeable.

In applications where the factor levels have a spatial or temporal structure, it is best to replace the equivalence matrix $E$ in (1) by a more suitable covariance matrix or generalized covariance function, justifying the neoclassical label. As an extreme example, consider a quantitative covariate, which is simply a factor taking values in the real line. The neoclassical Gaussian model with stationary additive random effects has the form

$$E(Y_i) = \beta_0 + \beta_1 x_i,$$
$$\text{cov}(Y_i, Y_j) = \sigma^2 \delta_{ij} + \sigma_s^2 K(|x_i - x_j|),$$

in which $K$ is a covariance function or generalized covariance function. Exchangeability implies $K(x, x') = \text{const} + \delta(x, x')$, but the more usual choices are Brownian motion in which $K(x, x') = -|x - x'|$, or integrated Brownian motion with $K(x, x') = |x - x'|^3$. The latter is a spline-smoothing model having the property that the predicted mean $E(Y(i^*)|\text{data})$ for a new unit such that $x(i^*) = x$ is a cubic spline in $x$ [Wahba (1990)]. This example may seem far removed from the sorts of factorial designs discussed in the paper, but factor levels are frequently ordered or partially ordered, in which case the argument for exchangeability of effects is not compelling. In principle, one may construct a similar covariance function for the effects of a conventional factor whose levels are ordered or partially ordered. Another option is to assume that the departures from linearity are exchangeable.

In time-series analysis, the spectrum determines a decomposition of the total sum of squares into components, two degrees of freedom for each Fourier frequency. Although there are no factors with identifiable levels in the conventional sense, by any reasonable interpretation of the term, this decomposition is an analysis of variance. In fact the key computational idea in the fast Fourier transform has its roots in Yates' algorithm for $2^n$ factorial designs, so the similarities are more than superficial. With this in mind, it is hard to understand Gelman's claim in Section 8 that analysis of variance is fundamentally about multilevel modeling. The canonical decomposition of the whole space as the direct sum of two-dimensional subspaces, one for each frequency, is a consequence of stationarity, or the group of translations. Any connection with exchangeability or the batching of coefficients is purely superficial.

Gelman's paper is a courageous attempt to reformulate a central part of applied statistics emphasizing Bayesian hierarchical modeling. Anyone who has taught factorial design and analysis at the graduate level will understand the constant difficult and sometimes painful struggle to achieve a reasonable and balanced attitude to the subject with its myriad and varied applications. Initially one tries to distill rules and extract common threads from typical applications, only to find later that all applications are atypical in one way or another. My own experience is that the state of this battle evolves as a process: it may converge, but it is not convergent to a fixed attitude or state of understanding. What seems important at one time often declines into insignificance later. It is clear that Gelman has thought hard about factorial models and ANOVA, and his views have evolved through consulting and teaching over a period of 10–15 years. My hope is that he will continue to think hard about the topic, and my prediction is that his views will continue to evolve for some time to come.

## REFERENCES

Cox, D. R. (1984). Interaction (with discussion). *Internat. Statist. Rev.* **52** 1–31. MR967201

Cox, D. R. and Snell, E. J. (1981). *Applied Statistics*: *Principles and Examples*. Chapman and Hall, London.

Joe, H. (1990). Extended use of paired comparison models with application to chess rankings. *Appl. Statist.* **39** 85–93. MR1038891

McCullagh, P. (2000). Invariance and factorial models (with discussion). *J. R. Stat. Soc. Ser. B Stat. Methodol.* **62** 209–256. MR1749537

Stewart, J. Q. (1948). Demographic gravitation: Evidence and application. *Sociometry* **11** 31–58.

Stigler, S. M. (1994). Citation patterns in the journals of statistics and probability. *Statist. Sci.* **9** 94–108.

Tukey, J. W. (1974). Named and faceless values: An initial exploration in memory of Prasanta C. Mahalanobis. *Sankhyā Ser. A* **36** 125–176.

Wahba, G. (1985). A comparison of GCV and GML for choosing the smoothing parameter in the generalized spline smoothing problem. *Ann. Statist.* **13** 1378–1402. MR811498

Wahba, G. (1990). *Spline Models for Observational Data*. MR1045442 SIAM, Philadelphia.

Wu, C. F. J. and Hamada, M. (2000). *Experiments*: *Planning, Analysis, and Parameter Design Optimization*. Wiley, New York. MR1780411

Department of Statistics
University of Chicago
5734 University Avenue
Chicago, Illinois 60637-1514
USA
e-mail: pmcc@galton.uchicago.edu